\documentstyle{article}
\begin{document}

\centerline {CONNECTION OF THE DIFFERENTIAL-GEOMETRICAL
STRUCTURES} \centerline {WITH SKEW-SYMMETRIC DIFFERENTIAL FORMS. }
\centerline {FORMING DIFFERENTIAL-GEOMETRICAL STRUCTURES AND
MANIFOLDS}
\bigskip
\centerline {L.I. Petrova}
\bigskip
\centerline{{\it Moscow State University, Russia, e-mail: ptr@cs.msu.su}}
\bigskip
The closure conditions of the inexact exterior differential form and 
dual form (an equality to zero of differentials of these forms) can be 
treated as a definition of some differential-geometrical structure. Such 
a connection discloses the properties and specific features of the 
differential-geometrical structures.

The using of skew-symmetric differential forms enables one to reveal a 
mechanism of forming differential-geometrical structures as well. To do 
this, it is necessary to consider the skew-symmetric differential forms 
whose basis, unlike to the case of exterior differential forms, are 
manifolds with unclosed metric forms. Such differential forms possess 
the evolutionary properties. They can generate closed differential forms, 
which correspond to the differential-geometrical structures. 
\bigskip

\section{Connection between exterior differential forms and
differential-geometrical structures} 

In this section it will be shown a connection between exterior differential 
forms and some differential-geometrical structures. This connection allows  
to describe the properties and specific features of differential-geometrical 
structures. 

The properties of differential-geometrical structures are based on 
the properties of closed exterior differential forms. 

\subsection*{Closed exterior differential forms} 

The exterior differential form of degree $p$ ($p$-form on the 
differentiable manifold) is called a closed one if its differential is 
equal to zero: 
$$
d\theta^p=0\eqno(1)
$$
(In more detail about skew-symmetric differential forms one can 
read in [1-3]. With the theory of exterior differential forms one 
can become familiar from the works [4-8]). 

From condition (1) one can see that the closed form is 
a conservative quantity. 

The differential of the form is a closed form. That is, 
$$
dd\omega=0\eqno(2)
$$
where $\omega$ is an arbitrary exterior form. 

The form, which is a differential of some other form: 
$$
\theta^p=d\theta^{p-1}\eqno(3)
$$
is called an {\it exact} form. Exact forms prove to be closed 
automatically [6] 
$$
d\theta^p=dd\theta^{p-1}=0\eqno(4)
$$

Here it is necessary to pay attention to the following points. In the 
above presented formulas it was implicitly assumed that the 
differential operator $d$ is a total one (that is, the 
operator $d$ acts everywhere in the vicinity of the point 
considered locally),  and therefore, it acts on the manifold of the 
initial dimension $n$. However, a differential may be internal. Such a 
differential acts on some structure with a dimension being less 
than that of the initial manifold. The structure, on which the exterior 
differential form may become a closed {\it inexact} form, is a 
pseudostructure with respect to its metric properties. \{Cohomology, 
sections of cotangent bundles, the eikonal surfaces, the 
characteristical and potential surfaces, 
and so on may be regarded as examples of pseudostructures.\} 

If the form is closed on  pseudostructure only, the closure 
condition is written as 
$$
d_\pi\theta^p=0\eqno(5)
$$
And the pseudostructure $\pi$ is defined from the condition 
$$
d_\pi{}^*\theta^p=0\eqno(6)
$$
where ${}^*\theta^p$ is the dual form. 
(For the properties of dual forms see [8]). 

From conditions (5) and (6) one can see that the form closed 
on pseudostructure (a closed inexact form) is a conservative object, 
namely, this quantity conserves on pseudostructure. 

\subsection*{Differential-geometrical structures} 

From the definition of a closed inexact exterior form one can see 
that to this form there correspond two conditions: 

1) condition (5) is a closure condition of the exterior form itself, 
and

2) condition (6) is that of the dual form. 

Conditions (5) and (6) can be regarded as a definition of some binary 
object that combines the pseudostructure and the exterior differential 
form being defined on this pseudostructure. Such a binary object 
can be named a Bi-Structure. The Bi-Structure combines both an algebraic 
object, namely, the closed exterior form, and a geometric one, namely, 
the pseudostructure as well. 

In its properties this differential-geometrical structure is a 
well-known G-Structure. Here a new term Bi-Structure 
was introduced to distinguish it from G-structures 
to which there correspond closed inexact exterior forms. 

The specific feature of this structure consists in the fact that 
it combines objects, which possess duality. The closed exterior 
differential form and the closed dual form are such objects. The 
existence of one object implies that the other one exists as well. 
It seems to be senseless to combine mathematically these two 
binary objects. However, this combined structure constitutes a 
unified whole that carries a double meaning. This statement will 
we explained later while considering a mechanism of forming the 
differential-geometrical structures. 

From conditions (5) and (6) one can see that Bi-Stricture 
constitute a conservative object, namely, a quantity that is 
conservative on the pseudostructure. 
Hence, such an object (Bi-Structure) can correspond to some 
conservation law. (It is from such object that the physical fields 
and corresponding manifolds are formed.) 

The connection of the differential-geometrical structures with exterior 
forms allows to understand the properties and specific features of the 
differential-geometrical structures.

\section{The properties of exterior differential forms and 
differential-geometrical structures} 

The basic properties of exterior differential forms are connected with 
the fact that any closed form is a differential. 

The exact form is, by definition, a differential (see condition (3)). 
In this case the differential is total. The closed inexact form is 
a differential too. The closed inexact form is an interior (on 
pseudostructure) differential, that is 
$$
\theta^p_\pi=d_\pi\theta^{p-1}\eqno(7)
$$

And so, any closed form is a differential of the form of a lower 
degree: the total one $\theta^p=d\theta^{p-1}$ if the form is exact, 
or the interior one $\theta^p=d_\pi\theta^{p-1}$ on pseudostructure if 
the form is inexact. (This may have the physical meaning: the form of 
lower degree can correspond to the potential, and the closed form by 
itself can correspond to the potential force.) 

\subsection*{Invariant properties of closed exterior differential forms 
and differential-geometrical structures. Nondegenerate transformations} 

Since the closed form is a differential, then it is evident 
that the closed form proves to be invariant under all 
transformations that conserve a differential. 

The examples of such nondegenerate transformations are unitary, tangent, 
canonical, and gradient transformations. 

To the nondegenerate transformations there are assigned closed forms 
of given degree. To the unitary transformations it is assigned (0-form), 
to the tangent and canonical transformations it is assigned (1-form), 
to the gradient transformations it is assigned (2-form) and so on. 
It should be noted that these transformations are {\it gauge 
transformations} for spinor, scalar, vector, tensor (3-form) fields. 
Hence one can see that the differential-geometrical Bi-Structure 
relate to the gauge type of the differential-geometrical G-Structure. 
They remain to be invariant under all transformations that conserve the 
differential. 

The closure of exterior differential forms, and hence their invariance, 
results from the conjugacy of elements of exterior or dual forms. 

From the definition of exterior differential form one can see that 
exterior differential forms have complex structure. Specific features 
of the exterior form structure are homogeneity with respect to the 
basis, skew-symmetry, integrating terms each including two objects of 
different nature (the algebraic nature for form coefficients, and the 
geometric nature for  base components). Besides,  the exterior form 
depends on the space dimension and on the manifold topology. The 
closure property of the exterior form means that any objects, namely, 
elements of the exterior form, components of elements, elements of 
the form differential, exterior and dual forms and others, turn out to 
be conjugated. (In the author's work [2] some types of conjugacy of 
exterior differential forms have been considered). A variety 
of objects of conjugacy leads to the fact that the closed forms can 
describe a great number of various differential-geometrical structures. 
In doing so, a connection of the  differential-geometrical structures 
with exterior forms enables one to understand what the invariant 
properties of the differential-geometrical structures are conditioned by. 

Since the conjugacy is a certain connection between two operators or 
mathematical objects, it is evident that, to express a conjugacy 
mathematically, it can be used relations. Just such relations constitute 
the basis of mathematical apparatus of the exterior differential forms. 
This is an identical relation. Identical relations of exterior 
differential forms disclose also the properties of 
differential-geometrical structures. 

\bigskip
\subsection*{Identical relations of exterior differential forms} 

Identical relations of exterior differential forms reflect the 
closure conditions of differential forms, namely, vanishing the form 
differential (see formulas (1), (5), (6)) and hence 
the conditions connecting the forms of consequent degrees (see formulas 
(3), (7)). Since the closure conditions of differential forms 
and dual forms specify the differential-geometrical structures, 
identical relations for exterior differential forms specify the 
differential-geometrical structures too. (Here it should be noted 
that formulas (1) and  (3) are not of interest for the 
differential-geometrical structures). 

The importance of the identical relations for exterior 
differential forms is manifested by the fact that practically in 
all branches of physics, mechanics, thermodynamics one faces such 
identical relations. The examples of such relations are the above 
presented Cauchy-Riemann conditions in the theory of complex 
variables, the transversality condition in the variational 
calculus, the canonical relations in  Hamilton formalism, formulas 
by Newton, Leibnitz, Green, the integral relations by Stokes, 
Gauss-Ostrogradskii, the thermodynamic and characteristical relations, 
the Bianchi identities, the eikonal relations and so on. (In 
more detail about identical (and nonidentical) relations it is 
outlined in the author's work [3]). 

The identical relations express the fact that each closed exterior 
form is a differential of some exterior form (with a degree less 
by one). In general form such an identical relation can be written as 
$$
d _{\pi}\phi=\theta _{\pi}^p\eqno(8)
$$
In this relation the form in the right-hand side has to be a 
{\it closed} one. (As it will be shown below, 
the identical relations are satisfied only on pseudostructures). 

In identical relation (8) in one side it stands the closed form and 
in other side does a differential of some differential 
form of the less by one degree, which is a closed form as well. 

The identical relations of another type are the analog of relation (8) 
obtained by differentiating or integrating this relation. 

Identical relations of exterior differential forms are a mathematical 
expression of various types of conjugacy that leads to closed exterior 
forms. They describe a conjugacy of any objects: the form elements, 
components of each element, exterior and dual forms, exterior forms of 
various degrees, and others. The identical relations, which are 
connected with different types of conjugacy, elucidate invariant, 
structural and group properties of exterior forms that are of great 
importance in applications. 

The functional significance of identical relations for exterior 
differential forms lies in the fact that they can describe a conjugacy 
of objects of different mathematical nature. This enables one to see 
internal connections between various branches of mathematics. Due to 
these possibilities the exterior differential forms, and correspondingly 
the differential-geometrical structures, have wide application in various 
branches of mathematics. 

\section{Mechanism of forming the differential-geometrical structures} 

It has been shown that the skew-symmetric closed exterior differential 
forms allow to describe the properties and specific features of the 
differential-geometrical structures. 

In this section it will be shown that the skew-symmetric 
differential forms describe also a process of forming the 
differential-geometrical structures. However, to do this, one must 
use skew-symmetric differential forms, which, in contrast to 
exterior (skew-symmetric) differential forms, possess the 
evolutionary properties, and for this reason they were named 
evolutionary differential forms. 

A peculiarity of the evolutionary differential forms consists in the fact 
that they generate exterior differential forms, which correspond to the 
differential-geometrical structures. This elucidates a process of 
forming the differential-geometrical structures. 

A radical distinction between the evolutionary forms and the exterior 
ones consists in the fact that the exterior differential forms are 
defined on manifolds with {\it closed metric forms}, whereas the 
evolutionary differential forms are defined on manifolds with {\it 
unclosed metric forms}.

\subsection*{Some properties of manifolds} 

In the definition of exterior differential forms a differentiable 
manifold was mentioned. 
But differentiable manifolds are not a single type of manifolds on 
which the exterior differential forms are defined. In the general 
case there are manifolds with structures of any types. The theory of 
exterior differential forms was developed just for such manifolds. 
They may be the Hausdorff manifolds, fiber spaces, the comological, 
characteristical, configuration manifolds and so on. Since all these 
manifolds possess structures of any types, they have one common property, 
namely, locally they admit one-to-one mapping into  Euclidean subspaces 
and into other manifolds or submanifolds of the same dimension [6]. 

When describing any processes in terms of differential equations, 
one has to deal with manifolds that do not allow one-to-one 
mapping described above. Such manifolds are, for example, manifolds 
formed by trajectories of elements of the system described by 
differential equations. The manifolds that can be called accompanying 
manifolds are variable deforming manifolds. Evolutionary 
differential forms can be defined on manifolds of this type. 

The closed metric forms define the manifold structure, and the commutators 
of metric forms define the manifold differential characteristics that 
specify the manifold deformation: bending, torsion, rotation, twist. 
The topological properties of manifolds are connected with the metric 
form commutators. The metric form commutators 
specify a manifold distortion. For example, the commutator of 
the zero degree metric form $\Gamma^\rho_{\mu\nu}$ characterizes the 
bend, that of the first degree form 
$(\Gamma^\rho_{\mu\nu}-\Gamma^\rho_{\nu\mu})$ characterizes the torsion, 
the commutator of the third-degree metric form $R^\mu_{\nu\rho\sigma}$ 
determines the curvature. 

It is evident that the manifolds, that are metric ones or possess the 
structure, have closed metric forms. It is with such manifolds that the 
exterior differential forms are connected. 

If the manifolds are deforming manifolds, this means that their 
metric form commutators are nonzero. That is, the metric forms of such 
manifolds turn out to be unclosed. The accompanying manifolds and 
manifolds appearing to be deforming ones are examples of such manifolds. 

The skew-symmetric evolutionary differential forms, whose basis are 
deforming manifolds, are defined on manifolds with unclosed metric forms. 

Thus, the exterior differential forms are skew-symmetric differential 
forms defined on manifolds, submanifolds or on structures with closed 
metric forms. Evolutionary differential forms are skew-symmetric 
differential forms defined on manifolds with metric forms that are 
unclosed. 

The evolutionary properties of the evolutionary skew-symmetric 
differential forms are just connected with properties of the metric 
form commutators. 

\subsection*{Specific features of the evolutionary differential form} 

The evolutionary differential form of degree $p$ ($p$-form) 
can be also written as an exterior differential form [2]. 

But the evolutionary form differential cannot be written similarly to 
that presented for exterior differential forms. In the evolutionary form 
differential there appears an additional term connected with 
the fact that the basis of the form changes [2]. 

For example, we again inspect the first-degree form 
$\omega=a_\alpha dx^\alpha$. 
[From here on the symbol $\sum$ will be omitted and it will be 
implied that a summation is performed over double indices.  Besides, the 
symbol of exterior multiplication will be also omitted for the 
sake of presentation convenience]. 

The differential of this form can 
be written as $d\omega=K_{\alpha\beta}dx^\alpha dx^\beta$, where 
$K_{\alpha\beta}=a_{\beta;\alpha}-a_{\alpha;\beta}$ are 
components of the commutator of the form $\omega$, and 
$a_{\beta;\alpha}$, $a_{\alpha;\beta}$ are the covariant 
derivatives. If we express the covariant derivatives in terms of 
the connectedness (if it is possible), they can be written 
as $a_{\beta;\alpha}=\partial a_\beta/\partial
x^\alpha+\Gamma^\sigma_{\beta\alpha}a_\sigma$, where the first 
term results from differentiating the form coefficients, and the 
second term results from differentiating the basis. (In 
Euclidean space covariant derivatives coincide with ordinary ones 
since in this case derivatives of the basis vanish). If 
we substitute the expressions for covariant derivatives into the 
formula for the commutator components, we obtain the following 
expression for the commutator components of the form $\omega$: 
$$
K_{\alpha\beta}=\left(\frac{\partial a_\beta}{\partial
x^\alpha}-\frac{\partial a_\alpha}{\partial
x^\beta}\right)+(\Gamma^\sigma_{\beta\alpha}-
\Gamma^\sigma_{\alpha\beta})a_\sigma\eqno(9)
$$
Here the expressions 
$(\Gamma^\sigma_{\beta\alpha}-\Gamma^\sigma_{\alpha\beta})$ 
entered into the second term are just the components of 
commutator of the first-degree metric form. 

The evolutionary form commutator of any degree contains the commutator 
of the manifold metric form of corresponding degree. The commutator of 
the exterior form does not contains a similar term because the 
commutator of metric form of manifold, on which the exterior form is 
defined, is equal to zero. 

\subsection*{Unclosure of evolutionary differential forms} 

The evolutionary differential form commutator, in contrast to that of 
the exterior one, cannot be equal to zero because it includes the metric 
form commutator being nonzero. This means that the evolutionary form 
differential is nonzero. Hence, the evolutionary differential form, in 
contrast to the case of the exterior form, cannot be closed. 

The commutators of evolutionary forms depend not only on the 
evolutionary form coefficients, but also on the characteristics of 
manifolds, on which this form is defined. As a result, such a dependence 
of the evolutionary form commutator produces the topological 
and evolutionary properties of both the commutator and the evolutionary 
form itself (this will be demonstrated below). 

Since the evolutionary differential forms are unclosed, the mathematical 
apparatus of evolutionary differential forms does not seem to possess 
any possibilities connected with invariant properties of closed exterior 
differential forms. However, the mathematical apparatus of evolutionary 
forms includes some new unconventional elements like 
nonidentical relations and degenerate transformations [3]. Just such 
peculiarities allow to describe evolutionary processes. 

\subsection*{Nonidentical relations of evolutionary differential forms 
and their specific features} 

The identical relations of closed exterior differential forms reflect 
a conjugacy of any objects. The evolutionary forms, being unclosed, 
cannot directly describe a conjugacy of any objects. But they allow a 
description of the process in which the conjugacy may appear (a 
process when closed exterior differential forms are generated). Such 
a process is described by nonidentical relations. 

Nonidentical relations  can be written as 
$$
d\psi \,=\,\omega^p \eqno(10)
$$
Here $\omega^p$ is the $p$-degree evolutionary form that is 
nonintegrable, $\psi$ is some form of degree $(p-1)$, and 
the differential $d\psi$ is a closed form of degree $p$. 

In the left-hand side of this relation it stands the form differential, 
i.e. a closed form that is an invariant object. In the right-hand 
side it stands the nonintegrable unclosed form that is not an invariant 
object. Such a relation cannot be identical. 

One can see a difference of relations for exterior forms and evolutionary 
ones. In the right-hand side of identical relation (see relation (8)) 
it stands a closed form, whereas the form in the right-hand side of 
nonidentical relation (10) is an unclosed one. 

Nonidentical relations are obtained while describing any processes. 
A relation of such type is obtained while analyzing the integrability 
of the partial differential equation. An equation is integrable 
if it can be reduced to the form $d\psi=dU$. However it 
appears that, if the equation is not subject to an additional 
condition (the integrability condition), it is reduced to the 
form (10), where $\omega$ is an unclosed form and it cannot be 
expressed as a differential. The first principle of thermodynamics is 
an example of nonidentical relation. 

While investigating real physical processes one often faces the 
relations that are nonidentical. But it is commonly believed that only 
identical relations can have any physical and mathematical meaning. 
For this reason one immediately attempts to impose 
a condition onto the nonidentical relation under which this relation 
becomes identical, and it is considered only a relation that can 
satisfy the additional conditions. And all remaining is rejected. 

In this approach the identical relations, which can be obtained from 
a given nonidentical relation, are obtained. So, the closed exterior 
forms and relevant differential-geometrical structures are found. 
However, this approach does not solve the evolutionary problem. In this 
approach there is no answer to the questions of how do the closed 
exterior differential forms appear, how do the structures originate. 
That is, in such approach there is no answer to the question of how 
the conditions, under which the closed exterior forms are obtained, 
are realized. It turns out that these conditions are realized while 
varying the nonidentical relation, which is a selfvarying relation. 

\subsection*{Selfvariation of the evolutionary nonidentical relation} 

The evolutionary nonidentical relation is selfvarying, 
because, firstly, it is nonidentical, namely, it contains 
two objects one of which appears to be unmeasurable, and, 
secondly, it is an evolutionary relation, namely, a variation of 
any object of the relation in some process leads to variation of 
another object and, in turn, a variation of the latter leads to 
variation of the former. Since one of the objects is an unmeasurable 
quantity, the other cannot be compared with the first one, and hence, 
the process of mutual variation cannot terminate. 

Varying the evolutionary form coefficients leads to varying the first 
term of the evolutionary form commutator (see (2.3)). In accordance with 
this variation it varies the second term, that is, the metric form of 
manifold varies. Since the metric form commutators specifies the 
manifold differential characteristics, which are connected with the 
manifold deformation (as it has been pointed out, the commutator of the 
zero degree metric form specifies the bend, that of second degree 
specifies various types of rotation, that of the third degree specifies 
the curvature), this points to the manifold deformation. This means that 
it varies the evolutionary form basis. In turn, this leads to variation 
of the evolutionary form, and the process of intervariation of the 
evolutionary form and the basis is repeated. Processes of variation of 
the evolutionary form and the basis are controlled by the evolutionary 
form commutator and it is realized according to the evolutionary relation. 

The process of the evolutionary relation selfvariation plays a governing 
role in description of the evolutionary processes. 

The significance of the evolutionary relation selfvariation consists in 
the fact that in such a process it can be realized conditions under 
which the identical relation is obtained from the nonidentical relation. 
These are conditions of degenerate transformation.

\subsection*{Degenerate transformations. Origination of 
differential-geometrical structures} 

To obtain the identical relation from the evolutionary nonidentical 
relation, it is necessary that a closed exterior differential form 
should be derived from the evolutionary differential form, which is 
included into evolutionary  relation. However, as it has been shown 
above, the evolutionary form cannot be a closed form. For this reason 
a transition from the evolutionary form is possible only to an 
{\it inexact} closed exterior form, which is defined on pseudostructure. 

To the pseudostructure it is assigned a closed dual form 
(whose differential vanishes). For this reason a transition 
from the evolutionary form to a closed inexact exterior form proceeds 
only when the conditions of vanishing the dual form differential are 
realized, in other words, when the metric form differential or 
commutator becomes equal to zero. 

Since the evolutionary form differential is nonzero, whereas the closed 
exterior form differential is zero, the transition from the evolutionary 
form to the closed exterior form is allowed only under {\it degenerate 
transformation}. The conditions of vanishing the dual form differential 
(the additional condition) are the conditions of degenerate 
transformation. 

Such conditions can just be realized under selfvariation of the 
nonidentical evolutionary relation. 

As it has been already mentioned, the evolutionary differential form 
$\omega^p$, involved into nonidentical relation (10) is an unclosed one. 
The commutator, and hence the differential, of this form is nonzero. 
That is,
$$
d\omega^p\ne 0\eqno(11)
$$
If the conditions of degenerate transformation are realized, then from 
the unclosed evolutionary form one can obtain a differential form closed 
on pseudostructure. The differential of this form equals zero. That is, 
it is realized the transition  

 $d\omega^p\ne 0 \to $ (degenerate transformation) $\to d_\pi \omega^p=0$,
$d_\pi{}^*\omega^p=0$

As conditions of degenerate transformation (additional conditions) 
it can serve any symmetries of the evolutionary form coefficients 
or its commutator. (While describing material system such additional 
conditions are related, for example, to degrees of freedom of the 
material system). 

Mathematically to the conditions of degenerate transformation there 
corresponds a requirement that some functional expressions become equal 
to zero. Such functional expressions are Jacobians, determinants, 
the Poisson brackets, residues, and others. 

Thus, while selfvariation of the evolutionary nonidentical 
relation the dual form commutator can vanish (due to the  symmetries of 
the evolutionary form coefficients or its commutator). This means that 
it is formed the pseudostructure on which the differential form turns 
out to be closed. The emergence of the form being closed on 
pseudostructure points out to origination of the 
differential-geometrical structures.

\subsection*{Obtaining identical relation from  nonidentical one} 

On the pseudostructure $\pi$ evolutionary relation (10) transforms into 
the relation 
$$
d_\pi\psi=\omega_\pi^p\eqno(12)
$$
which proves to be an identical relation. Indeed, since the form 
$\omega_\pi^p$ is a closed one, on the pseudostructure this form turns 
out to be a differential of some differential form. In other words, 
this form can be written as $\omega_\pi^p=d_\pi\theta$. Relation (12) 
is now written as 
$$
d_\pi\psi=d_\pi\theta
$$
There are differentials in the left-hand and right-hand sides of 
this relation. This means that the relation is an identical one. 

From evolutionary relation (10) it is obtained the identical on the 
pseudostructure relation. In this case the evolutionary relation itself 
remains to be nonidentical one. (At this point it should be 
emphasized that differential, which equals zero, is an interior one. 
The evolutionary form commutator becomes zero only on the 
pseudostructure. The total evolutionary form commutator is nonzero. That 
is, under degenerate transformation the evolutionary form differential 
vanishes only {\it on pseudostructure}. The total differential of the 
evolutionary form is nonzero. The evolutionary form remains to be 
unclosed.) 

It can be shown that all identical relations of the exterior 
differential form theory are obtained from nonidentical relations (that 
contain the evolutionary forms) by applying degenerate transformations. 

{\it The degenerate transform is realized as a transition to 
nonequivalent coordinate system: a transition from the accompanying 
noninertial coordinate system to the locally inertial that}. 
Evolutionary relation (10) and condition (11) relate 
to the system being tied to the accompanying manifold, whereas 
identical relations (12) may relate only to the locally inertial 
coordinate system being tied to a pseudostructure. 

Transition from nonidentical relation (10) to identical relation (12) 
means the following. Firstly, it is from such a relation that 
one can obtain the differential $d_\pi\psi$ and find the desired 
function $\psi_\pi$ (a potential). And, secondly, an emergence 
of the closed (on pseudostructure) inexact exterior form $\omega_\pi^p$ 
(right-hand side of relation (12)) points to an origination of the 
conservative object. This object is a conservative quantity (the closed 
exterior form  $\omega_\pi^p$) on the pseudostructure (the dual form 
$^*\omega^p$, which defines the pseudostructure). This object is an 
example of the differential-geometrical structure (Bi-Structure). 

This complex is a new conjugated object. Below it will be shown 
a relation between characteristics of these objects 
(the differential-geometrical structure) and characteristics of the 
evolutionary differential forms. 

Thus, the mathematical apparatus of evolutionary differential forms 
describes a process of generation of the closed inexact exterior 
differential forms, and this discloses a process of origination of 
the differential-geometrical structure, namely, a new conjugated object. 

It can be seen that the process of conjugating the objects and 
obtaining differential-geometrical structures is a mutual exchange 
between the quantities of different nature (for example, between 
the algebraic and geometric quantities, between the physical and 
spatial quantities) and vanishing some functional expressions 
(Jacobians, determinants and so on). This follows from the fact 
that a selfvariation of the nonidentical evolutionary relation and 
a transition from the nonidentical evolutionary relation  to 
identical one develop as a result of mutual variations of the 
evolutionary form coefficients (which have the algebraic nature) 
and the manifold characteristics (which have the geometric 
nature), and a realization of the degenerate transformation with 
obeying additional conditions. 

The evolutionary differential form is an unclosed form, that is, it is 
the form whose differential is not equal to zero. The differential of 
the exterior differential form equals zero. To the closed exterior form 
there correspond conjugated operators, whereas to the evolutionary form 
there correspond nonconjugated operators. A transition from the 
evolutionary form to the closed exterior form and origination of the 
differential-geometrical structures is a transition from nonconjugated 
operators to conjugated ones. This is expressed mathematically 
as a transition from a nonzero differential (the evolutionary form 
differential is nonzero) to a differential that equals zero (the closed 
exterior form differential equals zero).

\subsection*{Characteristics of Bi-Structure} 

Since the closed exterior differential form, which corresponds to the 
Bi-Structure arisen, was obtained from the nonidentical relation that 
involves the evolutionary form, it is evident that the Bi-Structure 
characteristics must be connected with those of the evolutionary form 
and of the manifold on which this form is defined, with the conditions 
of degenerate transformation and with the values of commutators of the 
evolutionary form and the manifold metric form. 

The conditions of degenerate transformation, as it was said before, 
determine the pseudostructures. The first term of the evolutionary form 
commutator determines the value of the discrete change (the quantum), 
which the quantity conserved on the pseudostructure undergoes when 
transition from one pseudostructure to another. The second term of the 
evolutionary form commutator specifies a characteristics that fixes the 
character of the initial manifold deformation, which took place before 
the Bi-Structure arose. (Spin is such an example). 

A discrete (quantum) change of a quantity proceeds in the direction 
that is normal (more exactly, transverse) to the pseudostructure. Jumps 
of the derivatives normal to the potential surfaces are examples of such 
changes. 

Bi-Structure may carry a physical meaning. Such binary objects are 
the physical structures from which the physical fields are formed. 
This has been shown by the  author in the works [1,9]. 

The connection of Bi-Structure with the skew-symmetric 
differential forms allows to introduce a classification of 
Bi-structure in dependence on parameters that specify the 
skew-symmetric differential forms and enter into nonidentical and 
identical relation of the skew-symmetric differential forms. To 
determine these parameters one has to consider the problem of 
integration of the nonidentical evolutionary relation. 

Under degenerate transformation from the nonidentical evolutionary 
relation one obtains a relation being identical on pseudostructure. 
Since the right-hand side of such a relation can be expressed in terms 
of differential (as well as the left-hand side), one obtains a relation 
that can be integrated, and as a result he obtains a relation with the 
differential forms of less by one degree. 

The relation obtained after integration proves to be nonidentical 
as well. 

The resulting nonidentical relation of degree $(p-1)$ (relation that 
includes the forms of the degree $(p-1)$) can be integrated once again 
if the corresponding degenerate transformation has been realized and 
the identical relation has been formed. 

By sequential integrating the evolutionary relation of degree $p$ (in 
the case of realization of the corresponding degenerate transformations 
and forming the identical relation), one can get closed (on the 
pseudostructure) exterior forms of degree $k$, where $k$ ranges 
from $p$ to $0$. 

In this case one can see that under such integration the closed (on the 
pseudostructure) exterior forms, which depend on two parameters, are 
obtained. These parameters are the degree of evolutionary form $p$ 
(in the evolutionary relation) and the degree of created closed 
forms $k$. 

In addition to these parameters, another parameter appears, namely, the 
dimension of space. If the evolutionary relation generates the closed 
forms of degrees $k=p$, $k=p-1$, \dots, $k=0$, to them there correspond 
the pseudostructures of dimensions $(N-k)$, where $N$ is the space 
dimension. \{It is known that to the closed exterior differential forms 
of degree $k$ there correspond skew-symmetric tensors of rank $k$ and to 
corresponding dual forms there do the pseudotensors of rank $(N-k)$, 
where $N$ is the space dimensionality. The pseudostructures correspond 
to such tensors, but only on the space formed.\}

\subsection*{The properties of pseudostructures and closed exterior 
forms. Forming fields and manifolds} 

As mentioned before, the additional conditions, namely, the conditions 
of degenerate transformation, specify the pseudostructure. But at every 
stage of the evolutionary process it is realized only one element of 
pseudostructure, namely, a certain minipseudostructure. The additional 
conditions determine a direction (a derivative of the function that 
specifies the pseudostructure) on which the evolutionary form 
differential vanishes. (However, in this case the total differential 
of the evolutionary form is nonzero). The closed exterior form is formed 
along this direction. 

While varying the evolutionary variable the minipseudostructures form 
the pseudostructure. The example of minipseudoctructure is the wave 
front. The wave front is an eikonal surface (the level surface), i.e. 
the surface with a conservative quantity. A direction that specifies the 
pseudostructure is a connection between the evolutionary and spatial 
variables. It gives the rate of changing the spatial variables. Such a 
rate is a velocity of the wave front translation. While its translation 
the wave front forms the pseudostructure. 

Manifolds with closed metric forms are formed by pseudostructures. They 
are obtained from manifolds with unclosed metric forms. In this case the 
initial manifold (on which the evolutionary form is defined) and the 
manifold with closed metric forms originated (on which the closed 
exterior form is defined) are different spatial objects. 

It takes place a transition from the initial manifold with unclosed 
metric form to the pseudostructure, namely, to the manifold with closed 
metric forms being created. Mathematically this transition 
(degenerate transformation) proceeds as {\it a transition from one 
frame of reference to another, nonequivalent, frame of reference.}

The pseudostructures, on which the closed {\it inexact} forms are 
defined, form the pseudomanifolds. (Integral surfaces, pseudo-Riemann 
and pseudo-Euclidean spaces are the examples of such manifolds). In this 
process dimensions of the  manifolds formed  are connected with the 
evolutionary form degree.

To transition from pseudomanifolds to metric manifolds it is assigned 
a transition from closed {\it inexact} differential forms to {\it exact} 
exterior differential forms. (Euclidean and Riemann spaces are examples 
of metric manifolds). 

Here it is to be noted that the examples of pseudometric spaces are 
potential surfaces (surfaces of a simple layer, a double layer and so 
on). In these cases the type of potential surfaces is connected with 
the above listed parameters. 

Since the closed metric form is dual with respect to some closed exterior 
differential form, the metric forms cannot become closed by themselves, 
independent of the exterior differential form. This proves that manifolds 
with closed metric forms are connected with the closed exterior 
differential forms. This indicates that the fields of conservative 
quantities are formed from closed exterior forms at the same time when 
the manifolds are created from the pseudoctructures. (The specific 
feature of the manifolds with closed metric forms that have been formed 
is that they can carry some information.) That is, the closed exterior 
differential forms and manifolds, on which they are defined, are 
mutually connected objects. On the one hand, this shows a duality of 
these two objects (the pseudostructure and the closed inexact 
exterior form), and, on the  other hand, this means that these objects 
constitute a unified whole. This whole is a new conjugated object 
(Bi-Structure).

1. Petrova L.~I., Exterior and evolutionary skew-symmetric differential 
forms and their role in mathematical physics. 

http://arXiv.org/pdf/math-ph/0310050

2. Petrova L.~I., Invariant and evolutionary properties of the 
skew-symmetric differential forms. 

http://arXiv.org/pdf/math.GM/0401039

3. Petrova L.~I., Identical and nonidentical relations. Nondegenerate 
and degenerate transformations. (Properties of the 
skew-symmetric differential forms).

http://arXiv.org/pdf/math.GM/0404109

4. Cartan E., Les Systemes Differentials Exterieus ef Leurs Application 
Geometriques. -Paris, Hermann, 1945. 

5. Bott R., Tu L.~W., Differential Forms in Algebraic Topology. 
Springer, NY, 1982. 

6. Schutz B.~F., Geometrical Methods of Mathematical Physics. 
Cambridge University Press, Cambridge, 1982. 

7. Encyclopedia of Mathematics. -Moscow, Sov.~Encyc., 1979 (in Russian). 

8. Wheeler J.~A., Neutrino, Gravitation and Geometry. Bologna, 1960. 

9. Petrova L.~I., Conservation laws. Their role in evolutionary 
processes. (The method of skew-symmetric differential forms). 

http://arXiv.org/pdf/math-ph/0311008

\end{document}